\documentclass[12pt]{amsart}
\usepackage{latexsym,amssymb,amsxtra,amsmath,amsfonts}
\usepackage{verbatim}
\usepackage{epsfig}
\usepackage[pdftex]{color} 
\usepackage[top=1in, bottom=1in, left=1in, right=1in]{geometry}
\usepackage[toc,page]{appendix}

\renewcommand{\theequation}{\arabic{equation}}

\newcommand{\beq}{\begin{equation}}
\newcommand{\eeq}{\end{equation}}
\newcommand{\beqa}{\begin{eqnarray}}
\newcommand{\eeqa}{\end{eqnarray}}
\newcommand{\beaa}{\begin{eqnarray*}}
\newcommand{\ben}{\begin{eqnarray*}}
\newcommand{\eaa}{\end{eqnarray*}}
\newcommand{\een}{\end{eqnarray*}}

\newcommand \nc {\newcommand}

\newtheorem{theorem}{Theorem}[section]
\newtheorem{lemma}[theorem]{Lemma}

\newtheorem{corollary}[theorem]{Corollary}

\nc \thref{Theorem \ref}
\nc \leref{Lemma \ref}
\nc \prref{Proposition \ref}
\nc \coref{Corollary \ref}
\nc \deref{Definition \ref}
\nc \exref{Example \ref}
\nc \reref{Remark \ref}

\newcommand\numberthis{\addtocounter{equation}{1}\tag{\theequation}}

\def\({\left(}
\def\){\right)}
\def\[{\left[}
\def\]{\right]}
\def\<{\left\langle}
\def\>{\right\rangle}

\begin{document}

\title[VOA recursion ]{
$\mathcal{W}$-algebra constraints and topological recursion for
$A_N$-singularity\\
[1ex]\footnotesize\mdseries
  (with an Appendix by Danilo Lewanski)
}

\author{Todor Milanov}

\address{T.M.: Kavli IPMU (WPI) \\ The University of Tokyo \\ Kashiwa \\ Chiba 277-8583 \\ Japan}
\email{todor.milanov@ipmu.jp}

\address{D.L.: Korteweg-de Vries Institute for Mathematics\\
University  of Amsterdam\\
Postbus 94248\\
1090 GE Amsterdam\\
The Netherlands} 
\email{D.Lewanski@uva.nl}

\thanks{{\em 2000 Math. Subj. Class.} 14D05, 14N35, 17B69}

\begin{abstract}
We derive a Bouchard--Eynard type topological recursion for the total
descendant potential of $A_N$-singularity. Our 
argument relies on a certain twisted representation of a Heisenberg
Vertex Operator Algebra (VOA) constructed via the periods of
$A_N$-singularity. In particular, our approach allows us to prove that
the topological recursion for the total descendant potential is
equivalent to a certain generating set of $\mathcal{W}$-algebra
constraints. 
\end{abstract}

\maketitle

\section{Introduction}

Motivated by his work in Gromov--Witten theory, Givental has
introduced the notion of a total descendant and a total ancestor
potential (see \cite{G1}). The definition makes sense for every
conformal semi-simple Frobenius manifold. The main input is the
so-called $R$-matrix and several copies of the Witten--Kontsevich
$\tau$-function normalised in an appropriate way (see
\cite{G1,G2}). On the other hand, it was proved by \cite{BOSS} and
\cite{M2} that the total ancestor potential can be reconstructed only
in terms of the $R$-matrix by using the {\em local} Eynard--Orantin
recursion. The main problem addressed in this paper is to find a
topological recursion for the total descendant potential. The first
step in solving this problem was suggested by Bouchard and Eynard in 
\cite{BE}. Their construction was successfully applied to obtain a recursion for
the total descendant potential of $A_N$-singularity in
\cite{DNOPS} (see Section 7). In general however, the method of Bouchard and Eynard is
not directly applicable, because the spectral curve is an infinite
sheet covering, i.e., not a Riemann surface. In this paper we would
like to suggest an approach based on the VOA construction of
\cite{BM}. We will focus on the case of $A_N$-singularity and hence we
will recover Theorem 7.3 in \cite{DNOPS}. Furthermore, our approach
allows us to compare the topological recursion and the
$\mathcal{W}$-constraints for the total descendant potential of
$A_N$-singularity (see \cite{BM}). More precisely, we prove that the so called {\em dilaton
 shift} identifies the differential operators of the topological
recursion with states in the $\mathcal{W}$-algebra corresponding to
the elementary symmetric polynomials. Constructing explicitly elements
of the $\mathcal{W}$-algebra is in general very difficult
problem. It would be interesting to find out other examples in which the
topological recursion can be used to construct generators of a
$\mathcal{W}$-algebra.

\subsection{Results}

Our main result will be stated entirely in terms of the root system of
type $A_N$. The formulation in terms of vertex algebras requires a
little bit more notation, so it will be given later on
in Section \ref{sec:3}. Let us fix the notation and recall the necessary background. 
Let $\mathfrak{h}\subset \mathbb{C}^{N+1}$ be the hyper-plane
$\chi_1+\cdots+\chi_{N+1}=0$, where $\chi_i$ are the standard
coordinate functions on $\mathbb{C}^{N+1}$ . Recall that the root
system of type $A_N$ can be realised as 
\ben
\Delta=\{\chi_i-\chi_j\ |\ 1\leq i\neq j\leq N\}\quad \subset \quad \mathfrak{h}^*.
\een
The corresponding Weyl group is the symmetric group on $N+1$ elements,
while its action on $\mathfrak{h}^*$ is induced from the standard
action on $(\mathbb{C}^{N+1})^*$ given by permuting
$\chi_1,\dots,\chi_{N+1}$.
Furthermore, the unique $W$-invariant bilinear form $(\ |\ )$ for
which $(\alpha|\alpha)=2$ for all $\alpha\in \Delta$ is induced  by
the following bilinear form on  $(\mathbb{C}^{N+1})^*$:
\ben
(\chi_i|\chi_j)=-\frac{1}{h}+\delta_{ij},
\een 
where $h:=N+1$. 

We define a set of differential operators on the infinitely many
variables 
\ben
\mathbf{t}=\{t_{k,a}\},\quad 1\leq a\leq N,\quad k\geq 0.
\een
Sometimes it is convenient to rescale the above variables and to work
with
\ben
x_{k,a}=\frac{ t_{k,a} } {(-a+h)(-a+2h)\cdots (-a+kh)},\quad 1\leq a\leq
N,\quad k\geq 0.
\een
First, we define a set of linear differential operators
\ben
\Phi_a(\lambda):= \sum_{m=0}^\infty \Big(
\lambda^m x_{m,a} \hbar^{-1/2} +
\lambda^{-m-1} (a+mh)
\hbar^{1/2} \partial/\partial x_{m,h-a}  \Big),\quad 1\leq a\leq h,
\een
where $\hbar$ is a formal parameter. Next we introduce the so called
{\em propagators}
\ben
P_{ij}(\lambda):=
\frac{
\eta^{i+j }
}{
(\eta^i -\eta^j )^2} \, \lambda^{-2},\quad  
1\leq i\neq j \leq h,
\een
where $\eta=e^{2\pi\sqrt{-1}/h}$. Finally, the differential operators
that we need are 
\ben
X_j(\lambda) = \sum_{a=1}^N
\eta^{-ja}\Phi_a(\lambda)\lambda^{-a/h},\quad 1\leq j\leq h
\een
and 
\beq\label{dop}
X_{j_1,\dots,j_r}(\lambda)  = \sum_{i_1,\dots,i_{r'}}
\Big(\prod_{s=1}^{r'} P_{i_s}(\lambda) \Big) \ 
{:} \prod_{j\in  J\setminus{I}} X_j(\lambda) {:},
\eeq
where the sum is over all disjoint pairs $i_s=(i_s^{(1)},i_s^{(2)})$,
$1\leq s\leq r'$, s.t., 
\ben
1\leq i_s^{(1)}<i_s^{(2)} \leq h, \quad
i_1^{(1)}<\cdots <i_{r'}^{(1)},
\een
we have used the notation
\ben
I=\bigcup_{s=1}^{r'} \{i_s^{(1)},i_s^{(2)}\},\quad
J=\{j_1,\dots,j_r\},\quad P_{i_s}(\lambda) = P_{ i_s^{(1)},i_s^{(2)} }(\lambda),
\een
and ${:}\, {:}$ is the normal ordering in which all differentiation
operations are applied before the multiplication ones.

The total descendant potential is a formal series of the type 
\ben
\mathcal{D}(\hbar;\mathbf{t}) = \exp\Big( \sum_{g=0}^\infty
\hbar^{g-1} \mathcal{F}^{(g)}(\mathbf{t}) \Big),
\een
where $\mathcal{F}^{(g)}$ are formal power series in $\mathbf{t}$. We
refer to \cite{G2} for the precise definition. Let us define
$\Omega^{(g)}_{j_1,\dots,j_r}$ by the following identity
\ben
X_{j_1,\dots,j_r} \, \mathcal{D}(\hbar;\mathbf{t}) = 
\Big( \sum_{g=0}^\infty \hbar^{g-r/2}
\Omega^{(g)}_{j_1,\dots,j_r}(\lambda;\mathbf{t}) 
\Big)
\mathcal{D}(\hbar;\mathbf{t}) ,
\een 
where $1\leq j_1<\cdots <j_r\leq h$.
\begin{theorem}\label{t1}
The following identity holds:
\ben
(-a+(m+1)h)
\frac{\partial \mathcal{F}^{(g)}}{\partial x_{m,a}} =
-\operatorname{Res}_{\lambda=0} \frac{1}{h}\sum_{i=1}^h \sum_{j_1,\dots,j_r}
\frac{
\eta^{-ia} \lambda^{m+1-\frac{1}{h}(a+r) } 
}{
\prod_{s=1}^r (\eta^i-\eta^{j_s})
}\, 
\Omega^{(g)}_{i,j_1,\dots,j_r}(\lambda;\mathbf{t})
d\lambda,
\een
where the 2nd sum is over all non-empty subsets $\{j_1,\dots,j_r\}$ of $\{1,\dots,i-1,i+1,\dots,h\}$. 
\end{theorem}
It is not hard to see that if we give an appropriate weight to each
variable $x_{k,i}$, so that the functions $\mathcal{F}^{(g)}$ are
homogeneous, then the identity in Theorem \ref{t1} will
give us a recursion that uniquely determines $\mathcal{F}^{(g)}$ for
all $g\geq 0$. 

\subsection{Genus-0}
Since the propagators do not contribute to genus 0, the genus-0
reduction of the identity in Theorem \ref{t1} takes a very simple
form. Put 
\ben
p_{m,a} = (-a+(m+1)h) \frac{\partial \mathcal{F}^{(0)}}{\partial
  x_{m,a}} ,\quad 1\leq a\leq N,\quad m\geq 0,
\een
\ben
\Phi^{(0)}_a(\lambda,\mathbf{t}):= \sum_{m=0}^\infty (x_{m,a}
\lambda^m + p_{m,h-a}\lambda^{-m-1}),
\een
and define the following numbers
\ben
C(a_1,\dots,a_r) = \sum_{1\leq j_1<\cdots <j_r\leq h-1} 
\frac{\eta^{-j_1a_1} }{1-\eta^{j_1}}
\dots
\frac{\eta^{-j_ra_r} }{1-\eta^{j_r }},\quad 1\leq a_1,\dots,a_r\leq N.
\een
\begin{corollary}\label{c1}
The following identity holds
\ben
p_{m,a} = -\operatorname{Res}_{\lambda=0}
\sum_{a_1,\dots,a_r=1}^{h-1}
C(a_1,\dots,a_r)\, \Phi^{(0)}_{a_0}(\lambda,\mathbf{t}) \Phi^{(0)}_{a_1}(\lambda,\mathbf{t})\cdots
\Phi^{(0)}_{a_r}(\lambda,\mathbf{t})
\lambda^{m+n+1}\, d\lambda,
\een
where the numbers $n\in \mathbb{Z}$ and $a_0,$ $0\leq a_0\leq h-1$ are
defined by
\ben
-(a+r+a_1+\cdots+a_r) = nh+a_0
\een
and if $a_0=0$ then we set $\Phi^{(0)}_{a_0}=0$.
\end{corollary}
If we set $x_{0,a}:=t_a$ and $x_{m,a}=0$ for $m>0$, then the identity
in Corollary \ref{c1} allows us to compute the primary potential of
the Frobenius structure. 

\subsection{$\mathcal{W}$-constraints}

Recall that the vector space
$\mathcal{F}:=\operatorname{Sym}(\mathfrak{h}[\zeta^{-1}]\zeta^{-1})$
has the structure of a highest weight $\widehat{\mathfrak{h}}$-module, where
$\widehat{\mathfrak{h}}:=\mathfrak{h}[\zeta,\zeta^{-1}] \oplus \mathbb{C}$ is
the Heisenberg Lie algebra with Lie bracket defined via the invariant
bi-linear form $(\ |\ )$ (see Section \ref{sec:3}). Following the
construction in \cite{BM} we define a state-field correspondence
$v\mapsto X(v)$, which to every $v\in \mathcal{F}$ associates a
{\em twisted field} $X(v)$. The latter is a differential operator on
a set of formal variables $q_{k,i}$, $1\leq i\leq N$, $k\geq 0$ whose
coefficients are Laurent polynomials in $\lambda^{1/h}$. Let us point
out that under the {\em dilaton shift}
\beq\label{d-shift}
t_{k,i} = q_{k,i}+\delta_{k,0}\delta_{i,N}, \quad 1\leq i\leq N,\quad k\geq 0,
\eeq  
the differential operators
\ben
X_{j_1,\dots,j_r}(\lambda) = X(\chi_{j_1}\cdots \chi_{j_r},\lambda),
\een
where $X(v,\lambda)$ denotes the value of $X(v)$ at the point
$\lambda$ and we identify $\mathfrak{h}\subset \mathcal{F}$ via $a\mapsto
a\,\zeta^{-1}.$   

Let $e_r\in \operatorname{Sym}(\mathfrak{h})$, $2\leq r\leq h$, be the
degree-$r$ elementary symmetric polynomials in
$\chi_1,\dots,\chi_h$. Note that from the topological recursion in
Theorem \ref{t1} we get a set of differential operators that
annihilates the total descendant potential
$\mathcal{D}(\hbar;\mathbf{t})$. 
\begin{theorem}\label{t2}
Under the dilaton shift \eqref{d-shift} the set of differential
constraints corresponding to the topological recursion turns into
\ben
\operatorname{Res}_{\lambda=0}\ \lambda^m X(e_{h+1-a},\lambda)\,
\mathcal{D}(\hbar;\mathbf{q}) = 0,\quad 1\leq a\leq N,\quad m\geq 0.
\een
\end{theorem}
The proof of Theorem \ref{t2} will be reduced to a combinatorial
identity, whose proof will be given in the Appendix. 
It is easy to check that all $e_r$, $2\leq r\leq h$, are in the kernel
of the screening operators $e^\beta_{(0)}$, $\beta\in \Delta$.
Therefore the main result in \cite{BM} and Theorem \ref{t2} give an
alternative proof of Theorem \ref{t1}. Let us point out that in
general the invariant polynomials are not in the
$\mathcal{W}$-algebra, so at least to the author, it is a little bit surprising that the
elementary symmetric polynomials have this property.

\section*{Acknowledgements} 
I am thankful to the organizers of the workshop ``Topological
Recursion and TQFTs'', MFO, Oberwolfach (2016), where part of this 
work was presented. I would like to thank also Bojko Bakalov, Guido
Carlet, and Sergey Shadrin for interesting discussions and especially
to Danilo Lewanski for several useful observations. 
This work is supported by JSPS Grant-In-Aid 26800003
and by the World Premier International Research Center Initiative (WPI Initiative),
MEXT, Japan. 

\section{Conformal Frobenius structure}
Let us recall the construction of a Frobenius structure on the space of
miniversal unfolding of $A_N$-singularity (see \cite{Du, He, SaT}).
Let 
\ben
F(s,x)=\frac{x^{N+1}}{N+1}+s_1 x^{N-1}+\cdots + s_N
\een
be a miniversal unfolding of singularity of type $A_N$. The
deformation parameters are allowed to take arbitrary complex values,
i.e., 
\ben
s=(s_1,\dots,s_N)\quad \in \quad B:=\mathbb{C}^N.
\een
The space $B$ is equipped with a semi-simple Frobenius structure as
follows. Using the so called {\em Kodaira--Spencer} isomorphism 
\beqa\label{KS-iso}
T_sB\cong \mathbb{C}[x]/(\partial_x F(s,x)),\quad \partial/\partial s_i
\mapsto \partial_{s_i} F\quad (\mbox{mod } \partial_x F)
\eeqa
we can equip each tangent space $T_sB$ with a multiplication $\bullet_s$ and with a
residue pairing
\beqa\label{res-pairing}
(\partial/\partial s_i,\partial/\partial s_j) =
\frac{1}{2\pi\sqrt{-1}} \oint_C \frac{ \partial_{s_i} F \partial_{s_j}
  F }{ \partial_x F} \, dx ,
\eeqa
where the contour of integration $C$ is a big loop enclosing the
critical points of $F$. The main property of the above pairing and
multiplication is that the family of connections
\beqa\label{DC}
\nabla = \nabla^{\rm LC} - z^{-1} \sum_{i=1}^N (\partial_{s_i}\bullet_s) ds_i
\eeqa
is flat. Here $z$ is a formal parameter, $\nabla^{\rm LC}$ is the Levi--Civita connection of the
residue pairing, and $\partial_{s_i}\bullet_s$ denotes the linear
operator in $T_sB$ of multiplication by the tangent vector
$\partial/\partial s_i$.

 The flatness of $\nabla$ implies the flatness of $\nabla^{\rm
   LC}$. We construct a trivialisation of the tangent and the
 cotangent bundle as follows. Let us denote by $H=\mathbb{C}[x]/x^N$
 the local algebra of $F(0,x)$. Then we have the following identifications
\ben
T^*B\cong TB \cong B\times T_0B\cong B\times H,
\een
where the first isomorphism is given by the residue pairing, the
second one uses the parallel transport with $\nabla^{\rm LC}$, and the
last one is the Kodaira--Spencer isomorphism. Let us choose a
flat coordinate system $t=(t_1,\dots,t_N)$, s.t., the point $t=0$
corresponds to $s=0$, and the vector fields
$\partial/\partial t_i$ correspond to the basis $\phi_i(x) = x^{N-i}$
$(1\leq i\leq N)$ of $H$.  

The connection \eqref{DC} can be extended also in the $z$-direction 
\beq\label{DC-z}
\nabla_{\partial/\partial z} = \frac{\partial}{\partial z} - \theta
z^{-1} + (E\bullet)z^{-2},
\eeq
where $\theta$ is the so called Hodge grading operator and $E$ is the
Euler vector field. Recall that via the Kodaira-Spencer isomorphism \eqref{KS-iso} $E$ corresponds to
$F$. In flat coordinates we have
\ben
E=\sum_{i=1}^N (1-d_i) t_i\partial/\partial {t_i},
\een
where $d_i= \operatorname{deg}(\phi_i) = (N-i)/(N+1)$ $(1\leq i\leq
N)$ is the so-called degree spectrum. The maximal degree
$D=d_1=(N-1)/(N+1)$ is called the {\em conformal dimension} of the
Frobenius manifold. The operator $\theta =
\frac{D}{2}-\operatorname{deg}$, i.e., 
\ben
\theta:H\to H,\quad \theta(\phi_i) = (D/2-d_i) \phi_i = \Big( -\frac{1}{2}+\frac{i}{N+1}\Big)\phi_i.
\een 
The connection operators \eqref{DC} and \eqref{DC-z} give rise to a
flat connection on the trivial bundle $(B\times \mathbb{C}^*)\times H
\to B\times \mathbb{C}^*,$ which is also known as the {\em Dubrovin's connection}.

\subsection{The periods of $A_N$-singularity}

Put $X=B\times \mathbb{C}$ and let 
\ben
\varphi: X\to B\times \mathbb{C},\quad \varphi(t,x) = (t, F(t,x)).
\een
The non-singular fibers $X_{t,\lambda}:=\varphi^{-1}(t,\lambda)$ form
a smooth fibration called the {\em Milnor fibration}. Let us denote by
$(B\times \mathbb{C})'$ the complement of the {\em discriminant }
locus of the map $\varphi$, i.e., $(B\times \mathbb{C})'$ is the set
of all $(t,\lambda)$, s.t., the fiber $X_{t,\lambda}$ consists of
$N+1$ pairwise distinct points. 

The periods of $A_N$-singularity are defined by 
\ben
I^{(n)}_a(t,\lambda) = -d_t \partial_\lambda^n \int_{a_{t,\lambda}}
d^{-1}\omega \quad \in T^*_tB\cong H,
\een
where $n$ is an arbitrary integer and $a\in
\mathfrak{h}:=\widetilde{H}_0(X_{0,1})\otimes \mathbb{C}$ is a reduced
homology cycle. The notation on the RHS is as follows: we denote by
$\omega=dx$ and $d^{-1}\omega=x$ (this is a $0$-form), the integration
cycle $a_{t,\lambda}$ is obtained from $a$ after choosing a reference
path in $(B\times \mathbb{C})'$ from $(0,1)$ to
$(t,\lambda)$ and using the parallel transport with respect to the
corresponding Gauss--Manin connection. Finally, $d_t$ is the De Rham differential on $B$. 

The period integrals are solutions to a connection $\nabla^{(n)} $,
which is a Laplace transform of the Dubrovin's connection
\ben
\nabla^{(n)}_{\partial_{t_i}} & = & \partial_{t_i} +
\frac{\phi_i\bullet}{\lambda-E\bullet}\,\Big(\theta-\frac{1}{2}-n\Big) ,\quad 1\leq i\leq N,\\
\nabla^{(n)}_{\partial_\lambda} & = & \partial_\lambda -
\frac{1}{\lambda-E\bullet}\,\Big(\theta-\frac{1}{2}-n\Big).
\een
The above system of equations can be solved in a neighbourhood of
$\lambda=\infty$ in the following way. Let us choose a solution to
Dubrovin's connection in the form $\Phi(t,z)=S(t,z) z^\theta$, where
$S(t,z)=1+S_1(t)z^{-1}+\cdots$ is an operator series whose
coefficients $S_k(t)\in \operatorname{End}(H)$. Such a solution is
unique and it satisfies the {\em symplectic condition}
$S(t,z)S(t,-z)^T=1$, where ${}^T$ is transposition with respect to the
residue pairing. The function 
\ben
Y^{(n)}(t,\lambda) = S(t,-\partial_\lambda^{-1})
\frac{\lambda^{\theta-n-1/2}}{\Gamma(\theta-n+1/2)} 
\een
is a fundamental solution to $\nabla^{(n)}$. Moreover, the reference
point $(0,1)$ is within the range of convergence (because $S(0,z)=1$),
therefore we can define an isomorphism $\mathfrak{h}\cong H$, s.t.,
\ben
I^{(n)}_a(t,\lambda) = Y^{(n)}(t,\lambda) a
\een
for all $(t,\lambda)$ sufficiently close to $(0,1)$.

\subsection{Monodromy representation}

Let us denote by $\Delta\subset H$ the set of vanishing cycles. 
\begin{lemma}\label{roots}
The set of vanishing cycles $\Delta=\{\chi_i-\chi_j\ |\ 1\leq i\neq
j\leq N+1\}$, where 
\ben
\chi_i=\sum_{a=1}^N \eta^{-ia}(N+1)^{-a/(N+1)}
\Gamma\Big(1-\frac{a}{N+1}\Big) \, \phi_{N+1-a},
\een
where $\eta=e^{2\pi\sqrt{-1}/(N+1)}.$
\end{lemma}
\proof
The fiber $X_{t,\lambda}$ consists of the zeroes $x_i(t,\lambda)$
$(1\leq i\leq N+1)$ of the equation $F(s,x)=\lambda$. The vanishing
cycles have the form $\alpha=[x_i(0,1)]-[x_j(0,1)]$, where $x_i(0,1) =
(N+1)^{1/(N+1)}\eta^i$. By definition 
\ben
I^{(0)}_\alpha(t,\lambda) = -d_t \int_{\alpha_{t,\lambda}} x = -d_t
(x_i(t,\lambda)-x_j(t,\lambda)).
\een
Furthermore,
\ben
-d_t\, x_i(t,\lambda) = \sum_{a=1}^N
\frac{x_i(t,\lambda)^{N-a}}{\partial_xF(t,x_i)} \, ds_a.
\een
On the other hand, note that the residue pairing has
the form
\ben
(\partial/\partial t_a,\partial/\partial t_b) = (x^{N-a},x^{N-b}) = \delta_{a+b,N+1}.
\een
Therefore, at $t=0$ we have 
\ben
I^{(0)}_\alpha(0,\lambda) =\sum_{a=1}^N (x_i(0,\lambda)^{-a}
-x_j(0,\lambda)^{-a} ) ds_a 
\een
and since at $t=0$: $ds_a=dt_a = x^{a-1}=\phi_{N+1-a}$, we get
\ben 
\sum_{a=1}^N (N+1)^{-a/(N+1)} (\eta^{-ia}-\eta^{-ja})
\lambda^{-a/(N+1)} \phi_{N+1-a} = Y^{(0)}(0,\lambda) (\chi_i-\chi_j).\qed
\een

Using the above Lemma we can verify Saito's formula for the
intersection pairing (see \cite{S0}), i.e., the bi-linear form 
\ben
(a|b):=(I^{(0)}_a(t,\lambda),(\lambda-E\bullet)I^{(0)}_b(t,\lambda))
\een
coincides with the intersection pairing in
$\widetilde{H}_0(X_{0,1})$. The Picard--Lefschetz formula for the
monodromy of the Gauss--Manin connection (see \cite{AGV}) takes the form
\ben
w_\alpha(y) = y-(\alpha|y)\alpha,\quad y\in H,
\een
where $w_\alpha$ is the image of the monodromy representation of
$\nabla^{(n)}$ 
\ben
\pi_1(B\times \mathbb{C})'\to \operatorname{GL}(H)
\een 
of a simple loop around the discriminant corresponding to a path along
which the cycle vanishes. In particular, we get that the monodromy
group is the symmetric group $S_{N+1}$ acting by permutation on the
set $(\chi_1,\dots,\chi_{N+1})$, while $w_\alpha$ for
$\alpha=\chi_i-\chi_j$ is just the transposition swapping $\chi_i$
and $\chi_j$.

\section{Heisenberg Vertex Operator Algebra}\label{sec:3}

Let us denote by $\widehat{\mathfrak{h}}$ the Heisenberg Lie algebra 
$H[\zeta,\zeta^{-1}]\oplus\mathbb{C}$ with bracket
\ben
[f(\zeta),g(\zeta)] = \operatorname{Res}_{\zeta=0}
(f'(\zeta)|g(\zeta)) d\zeta.
\een
It is convenient to denote $a_{(n)} = a\zeta^n$ for $a\in H$ and $n\in
\mathbb{Z}$. Then the above formula is equivalent to
\ben
[a_{(m)},b_{(n)}] = m(a|b)\delta_{m+n,0},\quad m,n \in \mathbb{Z}.
\een
The vector space $\mathcal{F} =
\operatorname{Sym}(H[\zeta^{-1}]\zeta^{-1})$ has a natural structure
of a highest-weight $\widehat{\mathfrak{h}}$-module, s.t., $a_{(n)}1=0$ for all $a\in
H$ and $n\geq 0$. 

\subsection{The tame Fock space}
Given a commutative ring $R$, let us denote by
$\widehat{\mathbb{V}}_R$ the space of formal series
of the form
\ben
\sum_{g\in\frac{1}{2}\mathbb{Z}} 
\sum_{K=((k_1,i_1),\dots,(k_s,i_s))}
 c^{(g)}_{K,I} \hbar^{g-1} t_{k_1,i_1}\cdots t_{k_s,i_s},\quad c^{(g)}_{K,I}\in R,
\een
where the 2nd sum is over all lexicographically increasing sequences
$K$ of pairs $(k,i)$, $k\geq 0, 1\leq i\leq N$, i.e., either $k_p<
k_{p+1}$ or $k_p=k_{p+1}$ and $i_p\leq i_{p+1}$
If $R=\mathbb{C}$, then we
simply put $\widehat{\mathbb{V}}:=\widehat{\mathbb{V}}_\mathbb{C}$ .
Let us denote by $\mathbb{V}_{\rm tame}\subset \mathbb{V}$ the subspace of formal series 
satisfying the {\em tameness condition}: if $c^{(g)}_{K,I}\neq 0$,
then 
\ben
k_1+\cdots +k_s\leq 3g-3+s.
\een
Let us denote by $\mathcal{O}$ the algebra of holomorphic functions on
the monodromy covering space of $(B\times\mathbb{C})'$. A {\em twisted
field} on $(B\times \mathbb{C})$ is a $\mathbb{C}$-linear map
$\mathbb{V}_{\rm tame} \to \widehat{\mathbb{V}}_\mathcal{O}$. The space
of all twisted fields will be denoted by
$\operatorname{Hom}_\mathbb{C}(\mathbb{V}_{\rm tame},
\widehat{\mathbb{V}}_\mathcal{O})$. 

\subsection{Twisted representation}

Following Givental \cite{G2}, we introduce the symplectic vector space
$\mathcal{H}=H(\!(z^{-1})\!)$ with the symplectic form
\ben
\Omega(f(z),g(z)) = \operatorname{Res}_{z=0} (f(-z),g(z)) dz.
\een 
Recall, also the following quantisation rules
\ben
(\phi_i z^k)\sphat = -\hbar^{1/2} \partial_{t_{k,i}},\quad 
(\phi^i(-z)^{-k-1} )\sphat =\hbar^{-1/2} \, t_{k,i},
\een
where $\phi^i:=\phi_{N+1-i}$ is the dual to $\phi_i$ with respect to
the residue pairing. These rules extend by linearity to define a
representation of the Poisson Lie algebra of linear and constant
functions on $\mathcal{H}$. We define a {\em State-Field} correspondence
\ben
X:\mathcal{F} \to \operatorname{Hom}_\mathbb{C}(\mathbb{V}_{\rm tame},
\widehat{\mathbb{V}}_\mathcal{O})
\een
as follows
\ben
X(a\zeta^{-1}) := \phi_a(t,\lambda) :=
(\phi_a(t,\lambda;z))\sphat,\quad a\in
H,\quad n\in \mathbb{Z}_{\geq 0},
\een
where
\ben
\phi_a(t,\lambda;z) = \sum_{n\in \mathbb{Z}} I^{(n+1)}_a(t,\lambda)\, (-z)^n.
\een
For the remaining states the definition is such that 
\ben
X_t(a_{(-n-1)} v,\lambda) = \operatorname{Res}_{\lambda'=\lambda}
\Big(
X_t(a,\lambda')X_t(v,\lambda)\frac{d\lambda'}{(\lambda'-\lambda)^{n+1}} \Big),
\een
where we denoted by $X_t(v,\lambda)$ the value of the field $X(v)$ at
a point $(t,\lambda)\in (B\times \mathbb{C})'$.

More explicitly, if $v=\alpha^1_{-k_1-1}\cdots \alpha^r_{-k_r-1} 1\in
\mathcal{F}$, then the field $X(v)$ can be computed explicitly in
terms of the generating fields $\phi_{\alpha^i}(t,\lambda)$ and the so
called {\em propagators} 
\ben
P_{\alpha,\beta}^{(k)}(t,\lambda)\in \mathcal{O},\quad \alpha,\beta\in H,\quad k\in
\mathbb{Z}_{\geq 0}
\een
defined by the Laurent series expansion
\ben
\Omega(\phi^+_\alpha(t,\lambda_1;z),\phi_\beta(t,\lambda_2;z)) =
\frac{(\alpha|\beta)}{(\lambda_1-\lambda_2)^2} +\sum_{k=0}^\infty
P^{(k)}_{\alpha,\beta}(t,\lambda_2)(\lambda_1-\lambda_2)^k. 
\een
The formula for the field $X(v)$ is
reminiscent of the Whick formula
\begin{equation}\label{VOA-operation}
X_t(v,\lambda) = \sum_J\Big(
\prod_{(i,j)\in J}
\partial_\lambda^{(k_j)}\, P_{\alpha^i,\alpha^j}^{(k_i)}(t,\lambda)\Big)
\; {:} \Big( \prod_{l\in J'} \partial_\lambda^{(k_l)} X_{s,\lambda}(\alpha^l)  \Big) {:} \,,
\end{equation}
where $\partial_\lambda^{(k)}:=\frac{\partial_\lambda^k}{k!}$ and the
sum is over all collections\/ $J$ of disjoint ordered pairs\/
$(i_1,j_1),$ $\dots,$ $(i_s,j_s)$ $\subset \{1,\dots,r\}$ such that\/
$i_1<\cdots<i_s$ and\/ $i_l<j_l$ for all $l$, and\/ 
$J' = \{1,\dots,r\} \setminus \{i_1,\dots,i_s,j_1,\dots,j_s\}$.

It is proved in \cite{M1} that the analytic continuation of the propagators is compatible with the
monodromy action on $\alpha$ and $\beta$. Moreover, we have the
following explicit formulas
\ben
\Omega(\phi^+_\alpha(t,\lambda_1;z),\mathbf{f}_\beta(t,\lambda_2;z)) =
\frac{1}{\lambda_1-\lambda_2}\, (I^{(0)}_\alpha(t,\lambda_1),(\lambda_2-E\bullet)I^{(0)}_\beta(t,\lambda_2))
\een
and 
\ben
P^{(0)}_{\alpha,\beta}(t,\lambda) = \frac{1}{2}(
(\lambda-E\bullet)I^{(1)}_\alpha(t,\lambda),I^{(1)}_\beta(t,\lambda)),
\een
where
\ben
\mathbf{f}_\beta(t,\lambda;z)=\sum_{n\in \mathbb{Z}}
I^{(n)}_\beta(t,\lambda)\, (-z)^n.
\een
The monodromy representation extends naturally to $\mathcal{F}$. It
follows from formula \eqref{VOA-operation} that the analytic
continuation of $X_t(v,\lambda)$ in $(t,\lambda)$ is compatible (or
equivalent) to the monodromy action on $v$.

\subsection{Global recursion}

If $(t,\lambda)\in (B\times \mathbb{C})'$ and $c^1,\dots, c^r \in H$,
then we define 
\ben
\Omega^{(g)}_{c^1\cdots c^r}(t,\lambda;\mathbf{t})\ \in \ \mathbb{C}[\![t_0,t_1,t_2,\dots]\!]
\een
by the following equation
\ben
X_t(c^1_{(-1)} \dots c^r_{(-1)} 1,\lambda) \,
  \mathcal{A}_t(\hbar;\mathbf{t}) =
\sum_{g=0}^\infty \hbar^{g-\frac{r}{2}} \, \Omega^{(g)}_{c^1\cdots
  c^r}(t,\lambda;\mathbf{t}) 
 \mathcal{A}_t(\hbar;\mathbf{t}) .
\een
Recall also, that the total ancestor potential has the form
\ben
\mathcal{A}_t(\hbar;\mathbf{t}) = \exp\Big(
\sum_{n,g=0}^\infty \frac{\hbar^{g-1}}{n!}\langle \mathbf{t}(\psi),\dots,\mathbf{t}(\psi)\rangle_{g,n}(t)\Big)
\een
where $\mathbf{t}(\psi)=\sum_{k=0}^\infty \sum_{a=1}^N t_{k,a}\phi_a
\psi^k$ and the correlator has the form
\ben
\langle \phi_{a_1}\psi^{k_1},\dots,\phi_{a_n}\psi^{k_n}\rangle_{g,n}(t) =
\int_{\overline{\mathcal{M}}_{g,n}} 
\Lambda^t_{g,n}(\phi_{a_1},\dots,\phi_{a_n})\psi_1^{k_1}\cdots \psi_n^{k_n},
\een
where $\Lambda_{g,n}^t:H^{\otimes n}\to
H^*(\overline{\mathcal{M}}_{g,n};\mathbb{C})$ is a certain
Cohomological Field Theory defined through the Frobenius structure. 
According to the main result in \cite{M2}, the total ancestor
potential is uniquely determined by the following recursion:
\begin{equation}
\label{local-rec}
\sum_{n=0}^\infty \frac{1}{n!}\langle \phi_a\,\psi^m,\mathbf{t},\dots,\mathbf{t}\rangle_{g,n+1}(t)= \frac{1}{4}\sum_{i=1}^N {\rm
  Res}_{\lambda=u_i} \,
\frac{(I_{\beta_i}^{(-m-1)}(t,\lambda),\phi_a)
}{(I^{(-1)}_{\beta_i}(t,\lambda),1)} \ \Omega^{(g)}_{\beta_i,\beta_i}(t,\lambda;\mathbf{t})\ d\lambda,
\end{equation}
where $u_i$ $(1\leq i\leq N)$ are the critical values of $F(t,x)$ and
$\beta_i$ is a cycle vanishing over $\lambda=u_i$. 

Let $C$ be a loop that encloses all critical values. Motivated by the
work of Bouchard and Eynard \cite{BE}, we would like to compare the
RHS of \eqref{local-rec} with the following integral
\beq\label{rec-integral}
-
\frac{1}{2\pi\sqrt{-1}} \oint_{C} 
\sum_{i=1}^{N+1} \sum_{J}
\frac{(I^{(-m-1)}_{\chi_i}(t,\lambda),\phi_a)}{\prod_{j\in J}
  (I^{(-1)}_{\chi_i-\chi_j}(t,\lambda),1)} 
\Omega^{(g)}_{\chi_i,\chi_{j_1},\dots,\chi_{j_r} }(t,\lambda;\mathbf{t})d\lambda,
\eeq
where the 2nd sum is over all non-empty subsets $J\subset \{1,\dots,
N+1\}\setminus{\{i\}}$ and $j_1,\dots,j_r$ are the elements of $J$. 

\begin{theorem}
The RHS of the local recursion \eqref{local-rec} coincides with the
integral \eqref{rec-integral}.
\end{theorem}
\proof
The integral \eqref{rec-integral} can be evaluated with the residue theorem. It is a
sum of the residues at the critical values. Let us verify that the
residue at $\lambda=u_1$ coincides with the corresponding residue in
the local recursion \eqref{local-rec}. Similar argument applies to the
remaining critical values. We may assume that
$\beta=\chi_1-\chi_2$ is the cycle vanishing at $\lambda=u_1$. 

The terms with $r=1$ contribute to the residue only if the set
$J\cup \{i\}$ contains $1$ or $2$, otherwise
$(\chi_i|\beta)=(\chi_j|\beta)=0$ and the entire expression is
analytic at $\lambda=u_1$. The 2 terms for which $i=1$, $J=\{2\}$ and
$i=2$, $J=\{1\}$ add up to 
\ben
-\operatorname{Res}_{\lambda=u_1}
\frac
{(I^{(-m-1)}_{\chi_1-\chi_2}(t,\lambda),\phi_a)}
{I^{(-1)}_{\chi_1-\chi_2}(t,\lambda),1)}\, 
\Omega^{(g)}_{\chi_1\chi_2}(t,\lambda;\mathbf{t})d\lambda.
\een
However, 
$
-\Omega^{(g)}_{\chi_1\chi_2} = \frac{1}{4}(\Omega^{(g)}_{\beta,\beta}-\Omega^{(g)}_{\chi_1+\chi_2,\chi_1+\chi_2})
$
and since $(\chi_1+\chi_2|\beta)=0$, the form
$\Omega^{(g)}_{\chi_1+\chi_2,\chi_1+\chi_2}$ is analytic at
$\lambda=u_1$, so it does not contribute to the residue. The above
residue coincides with the residue contribution at $\lambda=u_1$ of
\eqref{local-rec}. 

We claim that the terms for which the set $J\cup \{i\}$ contains
precisely one of the elements $1$ or $2$ cancel with the terms for
which $J\cup\{i\}$ contains both $1$ and $2$. To avoid cumbersome
notation put $\chi_i:=\chi_i\zeta^{-1}\in \mathcal{F}$ and 
\ben
X_I(t,\lambda):=X_t(\chi_{i_1}\cdots \chi_{i_s}1,\lambda)
\een
where $I=\{i_1,\dots,i_s\}\subset \{1,2,\dots, N+1\}$.
Let us compute
\beq\label{type-1}
-
\operatorname{Res}_{\lambda=u_1}\, 
\sum_{s=1}^2 \sum_{i=1}^{N+1} \sum_{J: (J\cup\{i\})\cap \{1,2\}=\{s\}} 
\frac{(I^{(-m-1)}_{\chi_i}(t,\lambda),\phi_a)}{\prod_{j\in J}
  (I^{(-1)}_{\chi_i-\chi_j}(t,\lambda),1)} 
X_{J\cup\{i\}}(t,\lambda)\mathcal{A}_t(\hbar;\mathbf{t}).
\eeq
We may replace $X_{J\cup\{i\}}(t,\lambda) \mathcal{A}_t$ by 
\beq\label{factoriz}
X_{ (J\cup\{i\}) \setminus{\{s\}} }(t,\lambda) (\phi^+_{\chi_s}(t,\lambda;z))\sphat
\mathcal{A}_t 
\eeq
because the remaning terms do not contribute to the residue. Note that
\ben
(\phi^+_{\chi_s}(t,\lambda;z))\sphat\
\mathcal{A}_t =
-\hbar^{1/2}
\sum_{g,n=0}^\infty 
\frac{\hbar^{g-1}}{n!}
\langle \phi^+_{\chi_s}(t,\lambda;\psi),\mathbf{t},\dots,\mathbf{t}\rangle_{g,n+1}.
\een
Recalling
the local recursion \eqref{local-rec} the expression \eqref{factoriz} is
transformed into 
\ben
\frac{1}{4}\,\hbar^{1/2}
\sum_{k=1}^N \operatorname{Res}_{\lambda'=u_k}\,
\frac
{
\Omega(
  \phi^+_{\chi_s}(t,\lambda;z),\mathbf{f}_{\beta_k}(t,\lambda';z))
}
{
(I^{(-1)}_{\beta_k}(t,\lambda'),1)
}
X_{ (J\cup\{i\}) \setminus{\{s\}} }(t,\lambda)
X_t(\beta_k^2,\lambda') \mathcal{A}_t(\hbar;\mathbf{t}).
\een
On the other hand
\ben
\Omega(
  \phi^+_{\chi_s}(t,\lambda;z),\mathbf{f}_{\beta_k}(t,\lambda';z))=
\frac{1}{\lambda-\lambda'}(I^{(0)}_{\chi_s}(t,\lambda),(\lambda'-E\bullet)I^{(0)}_{\beta_k}(t,\lambda'))=\frac{(\chi_s|\beta_k)}{\lambda-\lambda'} +\cdots,
\een
where the dots stand for a term analytic at $\lambda'=\lambda$. The sum \eqref{type-1} turns into
\ben
&&
-
\frac{1}{4}\,\hbar^{1/2}
\sum_{k=1}^N
\operatorname{Res}_{\lambda=u_1}\, 
\operatorname{Res}_{\lambda'=u_k}\,
\sum_{s=1}^2 \sum_{i=1}^{N+1} \sum_{J: (J\cup\{i\})\cap \{1,2\}=\{s\}}
\frac{1}{\lambda-\lambda'}(I^{(0)}_{\chi_s}(t,\lambda),(\lambda'-E\bullet)I^{(0)}_{\beta_k}(t,\lambda'))\,  \\
&&
\frac{(I^{(-m-1)}_{\chi_i}(t,\lambda),\phi_a)}{(I^{(-1)}_{\beta_k}(t,\lambda'),1)\prod_{j\in J}
  (I^{(-1)}_{\chi_i-\chi_j}(t,\lambda),1)} 
X_{ (J\cup\{i\}) \setminus{\{s\}} }(t,\lambda)
X_t(\beta_k^2,\lambda') \mathcal{A}_t(\hbar;\mathbf{t})d\lambda'd\lambda.
\een
Note that if we compute first the residue with respect to
$\lambda=u_1$ we would get $0$. Furthermore, the two residue
operations commute unless $k=1$. If $k=1$, then 
\ben
\operatorname{Res}_{\lambda=u_1}\, 
\operatorname{Res}_{\lambda'=u_1}\, = \operatorname{Res}_{\lambda'=u_1}\, 
\operatorname{Res}_{\lambda=\lambda'} .
\een
Recalling the definition of the State-Field correspondence we get 
\ben
&&
-\frac{1}{4}\,\hbar^{1/2}
\operatorname{Res}_{\lambda=u_1}\, 
\sum_{s=1}^2 \sum_{i=1}^{N+1} \sum_{J: (J\cup\{i\})\cap \{1,2\}=\{s\}}
(\chi_s|\beta_1)\times\\
&&
\frac{(I^{(-m-1)}_{\chi_i}(t,\lambda),\phi_a)}{(I^{(-1)}_{\beta_1}(t,\lambda),1)\prod_{j\in J}
  (I^{(-1)}_{\chi_i-\chi_j}(t,\lambda),1)} 
X_t(\chi_{j_1}\dots\chi_{j_r}\beta_1^2,\lambda)
\mathcal{A}_t(\hbar;\mathbf{t}) d\lambda,
\een
where $j_1,\dots,j_r$ are the elements of the set $J':= (J\cup\{i\})
\setminus{\{s\}}$. Just like before we can replace
$-\frac{1}{4}\beta_1^2$ with $\chi_1\chi_2$. Rearranging the sum so
that the summation over $J'$ is first we get
\ben
&&
\hbar^{1/2}
\operatorname{Res}_{\lambda=u_1}
\sum_{J'=\{j_1,\dots,j_r\}} 
\left(
\sum_{s=1}^2 \sum_{j'\in J'\cup\{s\}}
(\chi_s|\beta)\times\right. \\
&&
\left.
\frac{(I^{(-m-1)}_{\chi_{j'}}(t,\lambda),\phi_a)}{(I^{(-1)}_{\chi_1-\chi_2}(t,\lambda),1)\prod_{j\in
    J'\cup \{s\}\setminus{\{j'\}} }
  (I^{(-1)}_{\chi_{j'}-\chi_j }(t,\lambda),1)} \right)\, X_{J'\cup\{1,2\}}(t,\lambda)\mathcal{A}_t(\hbar;\mathbf{t}),
\een
where the outer sum is over all subsets $J'$ that do not contain $1$
and $2$. Note that the sum over $s$ and $j'$ in the brackets yields
\ben
\frac{(I^{(-m-1)}_{\chi_{1}}(t,\lambda),\phi_a)}{(I^{(-1)}_{\chi_1-\chi_2}(t,\lambda),1)\prod_{j\in J'}
  (I^{(-1)}_{\chi_1-\chi_j}(t,\lambda),1)}  -
\frac{(I^{(-m-1)}_{\chi_2}(t,\lambda),\phi_a)}{(I^{(-1)}_{\chi_1-\chi_2}(t,\lambda),1)\prod_{j\in J'}
  (I^{(-1)}_{\chi_2-\chi_j}(t,\lambda),1)} + 
\een
\ben
\sum_{j'\in J'}\left(
\frac{1}
{
(I^{(-1)}_{\chi_{j'}-\chi_1 }(t,\lambda),1)
}
-
\frac{1}
{
(I^{(-1)}_{\chi_{j'}-\chi_2 }(t,\lambda),1)
}
\right)\,
\frac
{
(I^{(-m-1)}_{\chi_{j'}}(t,\lambda),\phi_a)
}
{
(I^{(-1)}_{\chi_1-\chi_2}(t,\lambda),1) 
\prod_{j\in J'\setminus{\{j'\}}} 
(I^{(-1)}_{\chi_{j'}-\chi_j }(t,\lambda),1)
}.
\een
The above sum is precisely
\ben
\sum_{j'\in J'\cup\{1,2\}} 
\frac
{
(I^{(-m-1)}_{\chi_{j'}}(t,\lambda),\phi_a)
}
{
\prod_{j\in J'\cup\{1,2\}\setminus{\{j'\}}} 
(I^{(-1)}_{\chi_{j'}-\chi_j }(t,\lambda),1)
}.
\een
Note that the sum over $J$ and $i$ of the terms in \eqref{rec-integral} for which $J\cup\{i\}$
contains precisely one of the elements $1$ or $2$ coincides with
\eqref{type-1}. While the above argument shows that the sum
\eqref{type-1} cancels with the sum over $J$ and $i$ of the terms in
\eqref{rec-integral} for which $J\cup\{i\}$ contains both $1$ and
$2$. Therefore our claim follows and the proof of the Theorem is completed.  
\qed

Theorem \ref{t1} is an immediate corollary of the above theorem,
because in the case of $A_N$-singularity the restriction of the total
ancestor potential $\mathcal{A}_t(\hbar;\mathbf{t})$ to $t=0$
coincides with $\mathcal{D}(\hbar;\mathbf{t})$. 

\subsection{Example}
Let us use Corollary \ref{c1} to compute the primary genus-0 potential
of the $A_3$-singularity. Put $p_a:=p_{0,a}$, $x_{m,a}=0$ for
$m>0$, and $t_a:=x_{0,a}=t_{0,a}$. Note that 
\ben
\Phi^{(0)}_a = t_a + p_{4-a}\lambda^{-1}.
\een
The identities in Corollary \ref{c1} yield
\ben
p_3 & = &  t_1t_3+\frac{1}{2}t_2^2, \\
p_2 & = & 2 t_2 t_3- (C(1,1)+C(1,2)+C(2,1)) t_1^2t_2,\\
p_1 & = & -t_1p_3 + \frac{3}{2}\, t_3^2
-C(1,1)t_1^2t_3-(C(1,2)+C(2,1)+C(2,2))t_1t_2^2  \\
& & 
-(C(1,3)+C(3,1))t_1^2t_3-C(1,1,1)t_1^4.
\een
A straightforward computation gives
\ben
&& 
C(1,1)= C(2,2)=0,\\
&&
C(1,2)=1/2,\quad C(2,1)=1/2,\quad C(1,3)=(\eta-1)/2,\quad
C(3,1)=(-\eta-1)/2,\\
&&
C(1,1,1)=-1/4.
\een
Therefore
\ben
p_3  =  t_1t_3+\frac{1}{2}t_2^2, \quad
p_2  =  2 t_2 t_3- t_1^2t_2,\quad
p_1  =  -\frac{3}{2}t_1t_2^2 + \frac{3}{2}\, t_3^2+\frac{1}{4}t_1^4,
\een
i.e., 
\ben
\frac{\partial F}{\partial t_3} =  t_1t_3+\frac{1}{2}t_2^2,\quad
\frac{\partial F}{\partial t_2} = t_2 t_3- \frac{1}{2}t_1^2t_2,\quad
\frac{\partial F}{\partial t_1} =-\frac{1}{2}t_1t_2^2 + \frac{1}{2}\, t_3^2+\frac{1}{12}t_1^4,
\een
where $F$ is the restriction of $\mathcal{F}^{(0)}$ to $t_{0,a}=t_a$,
$t_{m,a}=0$ for $m>0$. Now it is easy to find that 
\ben
F(t_1,t_2,t_3) = \frac{1}{2}(t_1t_3^2+t_2^2t_3)-\frac{1}{4}t_1^2t_2^2+\frac{1}{60}t_1^5.
\een

\section{The topological recursion and $\mathcal{W}$-constraints}

The goal in this section is to prove Theorem \ref{t2}. Note that the
differential operators corresponding to the topological recursion have
the form
\beq\label{DO-1}
\sum_{r=0}^{h-1} \ 
\operatorname{Res}_{\lambda=0}\ 
\sum_{i=1}^h 
\sum_{\substack{1\leq j_1<\cdots <j_r\leq h\\  j_s\neq i}} 
\frac{
(I^{(-m-1)}_{\chi_i}(0,\lambda),\phi_a) 
}{
\prod_{s=1}^r (I^{(-1)}_{\chi_i-\chi_{j_s}}(0,\lambda),1) 
}
\, 
\hbar^{(r-1)/2}
\,
X_0(\chi_i\chi_{j_1}\cdots \chi_{j_r},\lambda).  
\eeq
Note that by definition if $r=0$, then the product over $s$ is $1$ and
the corresponding contribution to the sum is $\partial_{t_{m,a}}$. To
avoid cumbersome notation we set $X(v,\lambda):=X_0(v,\lambda)$. 
It is convenient to rewrite the above differential operator in terms
of the cycles 
\ben
\gamma_a := h^{-a/h}\Gamma(1-a/h) \, \phi_{h-a},\quad 1\leq a\leq N.
\een
Note that $\chi_i=\sum_{a=1}^N \eta^{-ia}\gamma_a$ and that 
\ben
I^{(0)}_{\gamma_a} (0,\lambda) = (h\lambda)^{-a/h}\,\phi_{h-a}.
\een
We get
\beq\label{kernel-den}
(I^{(-1)}_{\chi_i-\chi_j}(0,\lambda),1) = (\eta^i-\eta^j)\,
I^{(-1)}_{\gamma_N}(0,\lambda),1) =  (\eta^i-\eta^j)\, (h\lambda)^{1/h}
\eeq
and
\beq\label{kernel-num}
(I^{(-m-1)}_{\gamma_a}(0,\lambda),\phi_a) = 
\frac{(h\lambda)^{m+1-a/h}}
{
(-a+h)\cdots (-a+(m+1)h)}.
\eeq
The differential operator \eqref{DO-1} takes the form
\ben
\sum_{r=0}^{h-1}
\operatorname{Res}_{\lambda=0}\ d\lambda
\frac{
(I^{(-m-1)}_{\gamma_a}(0,\lambda),\phi_a) 
}{
 (I^{(-1)}_{\gamma_N}(0,\lambda),1)^r 
}\,
\hbar^{(r-1)/2} \times \\
\sum_{a_0,\dots,a_r=1}^N 
\sum_{i=1}^h 
\sum_{\substack{1\leq j_1<\cdots <j_r\leq h\\  j_s\neq i}} 
\,
\eta^{-i(r+a+a_0+a_1+\cdots + a_r)}
\Big(
\prod_{s=1}^r 
\frac{\eta^{-(j_s-i)a_s}}{ 1-\eta^{j_s-i} } \Big)
X(\gamma_{a_0}\gamma_{a_1}\cdots \gamma_{a_r},\lambda).
\een
Shifting the summation indexes $j_s\mapsto j_s+i$ and summing over $i$
we get
\beq\label{DO-2}
\sum_{r=0}^{h-1}
\operatorname{Res}_{\lambda=0}\ 
hd\lambda
\frac{
(I^{(-m-1)}_{\gamma_a}(0,\lambda),\phi_a) 
}{
 (I^{(-1)}_{\gamma_N}(0,\lambda),1)^r
}\,
\hbar^{(r-1)/2} 
\sum_{a_1,\dots,a_r=1}^N 
C(a_1,\dots,a_r)
X(\gamma_{a_0}\gamma_{a_1}\cdots \gamma_{a_r},\lambda),
\eeq
where $a_0$ is such that $0\leq a_0\leq h-1$, $r+a+a_0+\cdots +a_r
\equiv 0\ ({\rm mod}\, h)$, we assume that $\gamma_{a_0}=0$ if
$a_0=0$, and 
\ben
C(a_1,\dots,a_r):=
\sum_{1\leq j_1<\cdots <j_r\leq h-1 } 
\,
\prod_{s=1}^r 
\frac{\eta^{-j_sa_s}}{ 1-\eta^{j_s}  },
\een
where for $r=0$ the RHS is by definition 1. Since the differential operator
$X(\gamma_{a_0}\cdots\gamma_{a_r},\lambda)$ is invariant under the
permutations of $(a_0,\dots,a_r)$ we can arrange the 2nd sum in
\eqref{DO-2} to be over all increasing sequences $a_0\leq a_1\leq
\cdots \leq a_r$, i.e., 
\beq\label{BE-state}
\sum'_{1\leq a_0\leq \cdots \leq a_r\leq N} C[a_0,\dots,a_r] 
X(\gamma_{a_0}\gamma_{a_1}\cdots \gamma_{a_r},\lambda),
\eeq
where ${}'$ means that we allow only sequences $(a_0,\dots,a_r)$ that satisfy the
condition 
\ben
r+a+a_0+\cdots + a_r\equiv 0\ ({\rm mod}\, h)
\een
and the numbers $C[a_0,\dots,a_r]$ are defined as follows. If $r=0$,
then we put $C[a_0]:=1$. Otherwise,
\beq\label{coeff}
C[a_0,\dots,a_r] := \sum_{i=0}^r
\frac{1}{m_i}
\operatorname{SymC}(a_0,\dots,\widehat{a_i},\dots, a_r),
\eeq
where $m_i$ denotes the multiplicity of $a_i$ in the sequence
$(a_0,\dots,a_r)$ and $\operatorname{SymC}$ is the symmetrisation of
$C$
\ben
\operatorname{SymC}(b_1,\dots,b_r)  = 
\frac{1}{|\operatorname{Aut}(b_1,\dots,b_r)|} 
\sum_{\sigma\in S_r} C(a_{\sigma(1)},\dots,a_{\sigma(r)}). 
\een
Let us fix a summand in the sum \eqref{BE-state}. The corresponding
sequence has the form 
\ben
(a_0,a_1,\dots,a_r) = (b_1,\dots,b_{r'},N,\dots,N),\quad b_i<N, \quad 1\leq
i\leq r'.
\een
Put $m=r+1-r'$. Since the dilaton shift is equivalent to shifting
\ben
\gamma_a\mapsto \gamma_a + (I^{(-1)}_{\gamma_N}(0,\lambda),1)\,
\hbar^{-1/2} \delta_{a,N}.
\een
our summand is transformed into 
\ben
C[b_1,\dots,b_{r'},\underbrace{N,\dots,N}_{m}] 
\sum_{m'=0}^m \,
{m\choose m'}\,
X(b_1\cdots b_{r'}\underbrace{N\cdots N}_{m'},\lambda)\, 
 (I^{(-1)}_{\gamma_N}(0,\lambda),1)^{m-m'} \hbar^{-(m-m')/2}.
\een
The key step now is the following identity.
\begin{lemma}\label{remove-N}
The following identity holds
\ben
C[b_1,\dots,b_{r},\underbrace{N,\dots,N}_{m}] = (-1)^m 
{[\sum_{i=1}^r  b_i]_h \choose m}
C[b_1,\dots,b_r], 
\een
where $[b]_h$ denotes the remainder of $b$ modulo $h$. 
\end{lemma}
The proof will be given in the appendix. Using this Lemma we get
\ben
C[b_1,\dots,b_{r'},\underbrace{N,\dots,N}_{m}] 
{m\choose m'} = 
C[b_1,\dots,b_{r'},\underbrace{N,\dots,N}_{m'}] 
(-1)^{m-m'}{[\sum_{i=1}^{r'} b_i]_h -m' \choose m-m'} .
\een
Note that in particular, the multiplicity $m$ of $N$ in the sequence
$(a_0,\dots,a_r)$ does not exceed $[\sum_{i=1}^{r'} b_i]_h$. The sum
\eqref{BE-state} can be written as follows:
\ben
\sum_{r'=0}^{r+1}
\sum'_{1\leq b_1\leq \dots \leq b_{r'}<N} 
\sum_{m'=0}^{r+1-r'}
C[b_1,\dots,b_{r'},\underbrace{N,\dots,N}_{m'}] 
X(b_1\cdots b_{r'}\underbrace{N\cdots N}_{m'},\lambda)\times \\
\times 
(-1)^{r+1-r'-m'}
{[\sum_{i=1}^{r'} b_i]_h -m' \choose  r+1-r'-m'}
 (I^{(-1)}_{\gamma_N}(0,\lambda),1)^{r+1-r'-m'} 
\hbar^{-(r+1-r'-m')/2},
\een 
where the ${}'$ in the summation over $(b_1,\dots,b_{r'})$ means that 
\beq\label{arith-constr}
r'-1+b_1+\cdots + b_{r'} + a \equiv  0\ ({\rm mod}\, h).
\eeq
Substituting the above expression in \eqref{DO-2}, changing the summation
index $r$ via $s=r+1-r'-m'$, and changing the order of the summation
we get
\ben
\sum_{r'=0}^h 
\sum'_{1\leq b_1\leq \dots \leq b_{r'}<N} 
\sum_{m'=0}^{[\sum_{i=1}^{r'}b_i]_h}
\operatorname{Res}_{\lambda=0}\ 
hd\lambda
\frac{
(I^{(-m-1)}_{\gamma_a}(0,\lambda),\phi_a) 
}{
 (I^{(-1)}_{\gamma_N}(0,\lambda),1)^{r'+m'-1}
}\,
\hbar^{-1+(r'+m')/2} \times \\
\times 
C[b_1,\dots,b_{r'},\underbrace{N,\dots,N}_{m'}] 
X(b_1\cdots b_{r'}\underbrace{N\cdots N}_{m'},\lambda)\times \\
\times
\sum_{s=0}^{[\sum_{i=1}^{r'}b_i]_h-m'} 
(-1)^s{ [\sum_{i=1}^{r'}b_i]_h-m' \choose s}.
\een
Note that the sum over $s$ on the 3rd line of the above formula is $0$
unless $m'=[\sum b_i]_h$. Note also that $r'+m'\leq h$, otherwise the
2nd line of the formula vanishes. Recalling \eqref{arith-constr} we
get that $r'+m'=h+1-a$. Using formulas \eqref{kernel-den} and
\eqref{kernel-num} we get 
\ben
hd\lambda
\frac{
(I^{(-m-1)}_{\gamma_a}(0,\lambda),\phi_a) 
}{
 (I^{(-1)}_{\gamma_N}(0,\lambda),1)^{r'+m'-1}
}\,
\hbar^{-1+(r'+m')/2}  = 
\mbox{const} \ \hbar^{(h-a-1)/2} \ d\lambda\, \lambda^m ,
\een
where the value of the constant is not important. We get that up to a
constant the dilaton shift transforms the differential operator
\eqref{DO-1} into 
\ben
\operatorname{Res}_{\lambda=0} d\lambda \lambda^m\, 
\sum'_{1\leq b_1\leq \cdots \leq b_{h+1-a}\leq N} 
C[b_1,\dots,b_{h+1-a}]X(\gamma_{b_1}\cdots \gamma_{b_{h+1-a}},\lambda),
\een
where the ${}'$ indicates that the sum is over $(b_1,\dots,b_{h+1-a})$,
s.t., $\sum b_i \equiv 0 ({\rm mod}\ h)$. 
Put $r=h+1-a$. We claim that 
\ben
\sum'_{1\leq b_1\leq \cdots \leq b_r\leq N} 
C[b_1,\dots,b_r]\gamma_{b_1}\cdots \gamma_r
\een
coincides with the elementary symmetric polynomial in
$\chi_1,\dots,\chi_h$ of degree $r$. Similarly to what we did in the
beginning of this Section we can rewrite the above sum as 
\ben
\sum_{1\leq i_1<\cdots < i_r\leq h} \Big(
\sum_{s=1}^r 
\frac{
\eta^{i_s(r-1)} 
}{
\prod_{\substack{t=1 \\ t\neq s}}^r
    (\eta^{i_s}-\eta^{i_t})}
\Big)
\chi_{i_1}\dots \chi_{i_r}.
\een
The coefficient in front of $\chi_{i_1}\cdots \chi_{i_r}$ is 1,
because if we introduce the Vandermonde matrix
$A_{s,t}:=\eta^{(s-1)i_t}$, $1\leq s,t\leq r$, then 
the sum in the brackets can be interpreted as the quotient of the expansion of
$\operatorname{det}(A)$ with respect to the last row and 
\ben
\operatorname{det}(A) = \prod_{1\leq s<t\leq r} (\eta^{i_t}-\eta^{i_s}).
\een

\appendix

\section{Proof of Lemma \ref{remove-N}}

\medskip

\centerline{by D. Lewanski}

\medskip
Recall the definition of the numbers 
\ben
C[a_1,\dots,a_r],\quad 1\leq r\leq h,\quad 1\leq a_i\leq N
\een
given by formula \eqref{coeff}. It is convenient to extend the above
definition by setting $C[a_1,\dots,a_r]=0$ for $r>h$. 
\begin{lemma}
Let $r\geq 1$ and $1\leq a_1\leq\cdots \leq a_r\leq N-1$ be an arbitrary
sequence. The following identity holds:
\begin{align*}\label{SymC}
&\sum_{m=0}^\infty 
\operatorname{SymC}[a_1, \dots, a_r, \underbrace{N, \dots,
  N}_\text{m}] (1 - Y)^m =\\ \numberthis  =
&\frac{1-Y^h}{h(1-Y)} 
\sum_{k_1,\dots,k_r=0}^\infty 
\sum_{\substack{ I \subset \{1, \dots N\} \\  |I|=r  }} 
\Big(
\prod_{j= 1}^r \eta^{-i_j(a_j - k_j)} 
\Big)
Y^{\sum_{i=1}^r k_i},
\end{align*}
where the 2nd sum on the RHS is over all sequences $I=(i_1,\dots,i_r)$
of pairwise different numbers. 
\begin{proof}
Let us use the notation $I\subset \{1,2,\dots,N\}$ to denote that $I$ is a sequence
$(i_1,\dots,i_r)$ of pairwise distinct numbers, while $I\subset
(1,2,\dots,N)$ is a subsequence, i.e., a sequence of increasing
numbers $i_1<\cdots <i_r$. Recalling the definition of
$\operatorname{SymC}$ we get
\begin{align*}
&
\sum_{m=0}^\infty 
\operatorname{SymC}[a_1, \dots, a_r, \underbrace{N, \dots,
  N}_\text{m}] (1 - Y)^m = \\ 
&
= 
\sum_{m=0}^\infty
\sum_{\substack{I \subset \{1, \dots, N\} \\ |I | = r}} 
\prod_{s=1}^r \frac{\eta^{-i_s a_s}}{1 - \eta^{i_s}} 
\sum_{\substack{J \subset (1, \dots, N) \setminus I \\ |I | = m}} 
\prod_{t=1}^m \frac{\eta^{j_s}}{1 - \eta^{j_s}} (1 - Y)^m \\ 
&
=  
\sum_{m=0}^\infty \sum_{\substack{I \subset \{1, \dots, N\} \\ |I | = r}} 
\prod_{s=1}^r \frac{ - \zeta^{i_s (a_s + 1)}}{1 - \zeta^{i_s}} 
\sum_{\substack{J \subset (1, \dots, N) \setminus I \\ |J | = m}}  
\prod_{t=1}^m \frac{1}{1 - \zeta^{j_t}} (Y - 1)^m
\end{align*}
where $\zeta = \eta^{-1}$. Observe that for the function
$$
f_{I}(x) := 
\prod_{i \in \{1, \dots, N\} \setminus I} \!\!\!\!\!\!(x - \zeta^i) = 
\frac{x^h - 1}{x - 1} \prod_{i \in I} \frac{1}{(x - \zeta^i)}
$$
 we have 
$$
\frac{1}{m!} 
\frac{\partial_Y^m  f_{I}(Y)}{f_{I}(Y)} \Bigg{|}_{Y=1} = 
\sum_{\substack{J \subset (1, \dots, N) \setminus I \\ | J| = m}}  
\prod_{t=1}^m\frac{1}{1 - \zeta^{j_t}}
$$
contracting the Taylor expansion the initial term is:
$$
\sum_{\substack{I \subset \{1, \dots, N\} \\ |I | = r}} \prod_{s=1}^r \frac{ - \zeta^{i_s (a_s + 1)}}{1 - \zeta^{i_s}} \frac{f_I(Y)}{f_I(1)} =  \frac{Y^h - 1}{h(Y-1)} \sum_{\substack{I \subset \{1, \dots, N\} \\ |I | = r}} \prod_{s=1}^r\frac{ - \zeta^{i_s (a_s + 1)}}{Y - \zeta^{i_s}}
$$
Substituting back $\eta = \zeta^{-1}$ and expanding in geometric power series in the variables $Y \eta^{i_s}$ proves the lemma.

\end{proof}
\end{lemma}

The statement in Lemma \ref{remove-N} is equivalent to the following
identity.
\begin{lemma}
We have
$$
\sum_{m=0}^\infty 
C[a_1, \dots, a_r, \underbrace{N, \dots, N}_{m}] (1 - Y)^m =
Y^{[ \sum_{i=1}^r a_i ]_h} 
C[a_1, \dots, a_r],
$$
where $[a]_h$ denotes the remainder of $a$ modulo $h$.
\begin{proof}
By definition
\begin{align*}
&
\sum_{m=0}^\infty C[a_1, \dots, a_r, \underbrace{N, \dots, N}_\text{m}] (1 -Y)^m =\\
&
\quad (1-Y) \sum_{m=0}^\infty 
\operatorname{SymC}[a_1, \dots, a_r, \underbrace{N, \dots,
  N}_\text{m}] (1 - Y)^m + \\ 
&
\quad\sum_{i=1}^r 
\operatorname{SymC}[a_1, \dots, \hat{a}_i , \dots, a_r, \underbrace{N, \dots, N}_\text{m}] (1 - Y)^m.
\end{align*}
Let us substitute Equation (\ref{SymC} ) in the right hand side: the factor $(1-Y)^{-1}$ cancels out in the first summand, while in the $i$-th summand can be expanded as $\sum_{k_i=0}^{\infty}\eta^{-0(a_i - k_i)} Y^{k_i}$. Thus the first summand collects all the subsets of $\{0, \dots, N\}$ of cardinality $r$ not containing zero while the second summand collects all the subsets containing zero with the same cardinality $r$. Hence we get
\begin{align*}
\frac{(1-Y^h)}{h} \sum_{k_1,\dots,k_r=0}^\infty 
\sum_{\substack{  I \subset \{0, \dots, N \} \\ |I|=r }} 
\prod_{j=1}^r \eta^{-i_j(a_j - k_j)} Y^{\sum k_i}
\end{align*}
Now the set $\{0,1, \dots N\}$ is symmetric with respect to the shift
$i_j \mapsto i_j + 1$ simultaneously for all $j$. This implies
$\eta^{-\sum (a_j - k_j) } = 1$, hence $\sum k_i = [\sum a_i]_h + hl$,
for $l \in \mathbb{Z}_{\geq 0}$.
The initial term can now be expanded in powers of $Y$ as
$$ 
Y^{[\sum a_i]_h}(1 - Y^h) \frac{1}{h} \sum_{l=0} c_l
(Y^h)^l
$$
Since the expression is polynomial in $Y$, we should have $c_l = c_{l+1} = c$ for all indexes $l \geq 0$.
We showed: 
$$
\sum_{m=0}^\infty 
C[a_1, \dots, a_r, \underbrace{N, \dots, N}_\text{m}] (1 - Y)^m = \, 
Y^{[\sum a_i]_h} \frac{c}{h}
$$
Now evaluating at $Y=1$ gives $c/h = C[a_1, \dots, a_r]$ as desired. 

\end{proof}
\end{lemma}

\end{document}